\documentclass[amsfonts,11pt]{amsart}
\RequirePackage{amsmath,amssymb,amsthm, amscd, comment,mathtools}
\usepackage[margin=1.2in]{geometry}
\usepackage{amsmath}
\usepackage{amsthm}
\usepackage{amssymb}
\usepackage{amscd}
\usepackage{mathrsfs}
\usepackage{graphicx}
\usepackage{rotating, graphicx}
\usepackage[all,cmtip]{xy}
\usepackage{enumerate}
\usepackage{tabularx}
\usepackage{tikz-cd}
\usepackage[thinlines]{easytable}
\usepackage{multirow}
\usepackage{mathtools}
\usepackage{mathrsfs}
\usepackage{changepage} 
\usepackage{paralist}

\usepackage{kantlipsum}
\usepackage{mathabx}
\usepackage{enumitem}
\usepackage{times}
\usepackage{subfig}
\usepackage{color}
\usepackage[pagebackref,breaklinks,colorlinks,linkcolor=red,anchorcolor=red,citecolor=blue]{hyperref}

\calclayout

\newtheorem{proposition}[subsection]{Proposition}
\newtheorem{corollary}[subsection]{Corollary}
\newtheorem{theorem}[subsection]{Theorem}
\newtheorem{lemma}[subsection]{Lemma}
\theoremstyle{definition}
\newtheorem{definition}[subsection]{Definition}
\theoremstyle{remark}

\newtheorem{remark}[subsection]{Remark}

\numberwithin{equation}{subsection}

\DeclareMathAlphabet{\mathbbold}{U}{bbold}{m}{n}

\title{The quadratic Artin conductor of a motivic spectrum}

\author{Fangzhou Jin}
\address{School of Mathematical Sciences\\
Tongji University\\
Siping Road 1239\\
200092 Shanghai\\
China}
\email{\href{mailto:fangzhoujin@tongji.edu.cn}{fangzhoujin@tongji.edu.cn}}
\urladdr{\url{https://fangzhoujin.github.io/}}

\author{Enlin Yang}
\address{School of Mathematical Sciences\\
Peking University\\
No.5 Yiheyuan Road Haidian District,\\
Beijing 100871,\\
P.R.China}
\email{\href{mailto: yangenlin@math.pku.edu.cn}{yangenlin@math.pku.edu.cn}}
\urladdr{\url{https://www.math.pku.edu.cn/teachers/yangenlin/ely.htm}}

\date{\number\day-\number\month-\number\year}

\begin{document}

\maketitle

\begin{abstract}

Given a constructible motivic spectrum over a smooth proper scheme which is dualizable over an open subscheme, we define its quadratic Artin conductor under some assumptions, and prove a formula relating the quadratic Euler-Poincar\'e characteristic, the rank and the quadratic Artin conductor. As a consequence, we obtain a quadratic refinement of the classical Grothendieck-Ogg-Shafarevich formula.

\end{abstract}

\setcounter{tocdepth}{1}
\tableofcontents

\noindent

\section{Introduction}

\subsection{}
Let $k$ be an algebraically closed field and let $C$ be a projective smooth curve over $k$. Let $\mathcal{F}$ be an \'etale sheaf over $C$ which is locally constant outside a finite number of points $x_1,\cdots, x_n$. The famous \emph{Grothendieck-Ogg-Shafarevich formula} (see \cite[X Thm. 7.1]{SGA5}, where it is called the \emph{Euler-Poincar\'e formula}), takes form
\begin{align}
\label{eq:EP}
\chi(\mathcal{F})=\operatorname{rk}(\mathcal{F})\cdot \chi(C)-\sum_{i=1}^nArt_{x_i}(\mathcal{F}),
\end{align}
where $Art_{x_i}(\mathcal{F})$ is the \emph{Artin conductor} of $\mathcal{F}$ at $x_i$; in particular, the Euler-Poincar\'e characteristic of $\mathcal{F}$ depends not only on the ranks of the stalks of $\mathcal{F}$, but also on the information on the wild ramification of $\mathcal{F}$ at the points $x_1,\cdots, x_n$. The formula~\eqref{eq:EP} is a pioneering step towards the study of ramifications of \'etale sheaves, and has been generalized to the higher-dimensional case (\cite[Thm. 4.2.9]{KS}, \cite[Thm. 7.13]{Sai}).

\subsection{}
The main purpose of this paper is to prove a quadratic refinement of the Grothendieck-Ogg-Shafarevich formula using motivic homotopy theory:
\begin{theorem}[Corollary~\ref{cor:GOS2} and Corollary~\ref{cor:GOS3}]
Let $k$ be a field and let $p:X\to \operatorname{Spec}(k)$ be a smooth and proper morphism with $X$ connected, and let $Z$ be a nowhere dense closed subscheme of $X$ with open complement $U$. Let $K\in\mathbf{SH}_c(X)$ be a constructible motivic spectrum over $X$ such that $K_{|U}$ is dualizable. Then the following equality holds in the 2-inverted Grothendieck-Witt group of quadratic forms over $k$:
\begin{align}
\label{eq:gosintro}
\chi(p_*K)
=
p_*(\operatorname{rk}(K)\cdot e(T_{X/k}))-Art(K)\in GW(k)[1/2].
\end{align}
If $k$ has characteristic different from $2$ and $X$ is odd-dimensional, then one has
\begin{align}
\label{eq:GOS3}
\chi(p_*K)
=
\operatorname{rk}(K_{\textrm{et}})\cdot \chi(X/k)-Art(K)\in GW(k).
\end{align}
\end{theorem}
Here the two terms $\chi(p_*K)$ and $p_*(\operatorname{rk}(K)\cdot e(T_{X/k}))$ are analogs of $\chi(\mathcal{F})$ and $\operatorname{rk}(\mathcal{F})\cdot \chi(C)$ above, respectively (see~\ref{num:goss}); the element $\operatorname{rk}(K_{\textrm{et}})\in\mathbb{Z}$ is the rank of the \'etale realization of $K$. The term $Art(K)$ is called the \textbf{quadratic Artin conductor} of $K$ (Definition~\ref{def:quadSwan}), whose definition is more involved and will be explained below. We expect that the formula~\ref{eq:gosintro} also holds in $GW(k)$ when $k$ has characteristic different from $2$ and $Z$ has everywhere codimension at least $2$ in $X$, but at the moment we do not know how to deal with the case where $X$ is even dimensional.

\subsection{}
Our main tool is the use of the Milnor-Witt spectrum $\mathbf{H}_{\rm{MW}}\mathbb{Z}$ (\cite[Ch. 6 Def. 4.1.1]{BCD+}). Compared to the sphere spectrum $\mathbbold{1}_k$, the Milnor-Witt spectrum has the advantage of being more computable and flexible while still providing much quadratic information. One of the first obstructions in obtaining a formula of the form~\eqref{eq:gosintro} lies in the definition of the \emph{rank} $\operatorname{rk}(K)$, which is no longer an integer as in the \'etale case: since $K_{|U}$ is dualizable, one may take the categorical Euler characteristic of $K_{|U}$, but such a construction only gives an endomorphism of $\mathbbold{1}_U$ instead of an endomorphism of $\mathbbold{1}_X$. We consider the $\mathbb{E}$-valued rank $\operatorname{rk}(K,\mathbb{E})$ (Definition~\ref{def:rankE}) for $\mathbb{E}=\mathbf{H}_{\rm{MW}}\mathbb{Z}$ or some of its variants (\ref{num:spE}), and use some vanishing properties (\S~\ref{sec:van}) to provide cases when it can be lifted to $X$, see~\ref{num:rkwd} for a summary.

\subsection{}
Let $f:X\to \operatorname{Spec}(k)$ be the structure morphism with $\delta:X\to X\times_kX$. Recall that in \cite[Def. 5.1.3]{JY}, we constructed a map
\begin{align}
\begin{split}
\label{eq:CX}
C_X(K,\mathbb{E}):\mathbbold{1}_X
&\to
\underline{Hom}(K,K)
\simeq
\delta^!(\mathbb{D}_{X/k}(K)\boxtimes_kK)
\to
\delta^*(\mathbb{D}_{X/k}(K)\boxtimes_kK)\\
&=
\mathbb{D}_{X/k}(K)\otimes_kK
\simeq
K\otimes_k\mathbb{D}_{X/k}(K)
\to
f^!\mathbbold{1}_k
\to
f^!\mathbb{E}
\end{split}
\end{align}
called the $\mathbb{E}$-valued \emph{characteristic class} of $K$. Our definition of the quadratic Artin conductor is based on a localized variant of the class $C_X(K,\mathbb{E})$, see \cite[\S5.2]{AS}. Using the vanishing properties discussed in \S~\ref{sec:van}, we show in Proposition~\ref{prop:locex} that the \emph{($\mathbf{H}_{\rm{MW}}\Lambda$-valued) localized characteristic class} 
$C_X^Z(K,\mathbf{H}_{\rm{MW}}\Lambda))\in\mathbf{H}_{\rm{MW}}\Lambda_0(Z/k)$, is well-defined (in the sense of Definition~\ref{def:swan}) under one of the following assumptions:
\begin{enumerate}
\item
$2$ invertible in $\Lambda$;
\item
$Z$ everywhere of codimension at least $2$ in $X$;
\item
$X$ is odd-dimensional.
\end{enumerate}
Here for each case we use different arguments: the first case follows from the properties of real schemes, the second uses vanishings of Milnor-Witt cohomology, and the third one is a consequence of the vanishing of the Witt-valued Euler class for bundles of odd rank (\cite[Lemma 4.3]{LevWitt}, \cite[Prop. 4.4]{LevA}).

\subsection{}
In Theorem~\ref{th:GOS}, we prove a Grothendieck-Ogg-Shafarevich-type formula for the class $C_X^Z(K,\mathbb{E})$:
\begin{align}
\label{eq:GOSCintro}
C_X(K,\mathbb{E})=\operatorname{rk}(K,\mathbb{E})\cdot C_X(\mathbbold{1}_X,\mathbb{E})+i_*C_X^Z(K,\mathbb{E})\in \mathbb{E}_0(X/k)
\end{align}
under very similar assumptions. The proof uses the arguments as in Proposition~\ref{prop:locex}, together with the case $\mathbb{E}=\mathbf{H}\Lambda$ the motivic Eilenberg-Mac Lane spectrum, for which we use the idea of the proof of \cite[Prop. 5.2.3]{AS}.

\subsection{}
Since $X$, and therefore $Z$, is assumed to be proper over $k$, we define $Art(K)$ as the degree of the class $C_X^Z(K,\mathbb{E})$. Taking the degree of the formula~\eqref{eq:GOSCintro} gives an equality in $GW(k)[1/2p]$, where $p$ is the exponential characteristic of $k$ (Corollary~\ref{cor:GOS2p}). If $k$ has characteristic different from $2$, we use an argument of Levine (\cite[Rem. 2.1]{LevA}) to get rid of $p$ in the denominator and obtain the formula~\eqref{eq:gosintro}.

\subsection{}
In order to study the functorialities of the localized characteristic class, we introduce a relative version in Section~\ref{num:NAclass}, and show that it satisfies desired base change, push-forward and pull-back functorialities.

\subsubsection*{\bf Acknowledgments}
We would like to thank Marc Levine, Peng Sun and Nanjun Yang for very helpful discussions.
Both authors are supported by the National Key Research and Development Program of China Grant Nr.2021YFA1001400. F. Jin is supported by the National Natural Science Foundation of China Grant Nr.12101455 and the Fundamental Research Funds for the Central Universities. E. Yang is supported by Peking University's Starting Grant Nr.7101302006 and NSFC Grants Nr.12271006 and Nr.11901008.

\section{Notations and preliminaries}

\subsection{}
If $f:X\to Y$ is a separated morphism of finite type and $\mathbb{E}\in\mathbf{SH}(Y)$, we denote $\mathbb{E}_X=f^*\mathbb{E}$ and $\mathbb{E}^!_{X/Y}=f^!\mathbb{E}$. We denote $\mathcal{K}_{X/k}=f^!\mathbbold{1}_k$.

\subsection{}
\label{num:spE}
Throughout the paper, we will denote by $\mathbb{E}$ one of the following motivic spectra: 
\begin{itemize}
\item
the Milnor-Witt spectrum $\mathbf{H}_{\rm{MW}}\Lambda$, 
\item
the Witt spectrum $\mathbf{H}_{\rm{W}}\Lambda$, 
\item
the real \'etale spectrum $\mathbf{H}_{\rm{ret}}\mathbb{Z}$,
\item
the motivic Eilenberg-Mac Lane spectrum $\mathbf{H}\Lambda$.
\end{itemize}
All these spectra are commutative ring spectra endowed with a unit section $\mathbbold{1}\to\mathbb{E}$.

%

\subsection{}
By \cite[Ch. 3 Thm. 5.0.2]{BCD+}, we have canonical isomorphisms
\begin{align}
\label{eq:MW2split}
\mathbf{H}_{\rm{MW}}\mathbb{Z}[1/2]
\simeq
\mathbf{H}\mathbb{Z}[1/2]
\oplus
\mathbf{H}_{\rm{W}}\mathbb{Z}[1/2]
\simeq
\mathbf{H}\mathbb{Z}[1/2]
\oplus
\mathbf{H}_{\rm{ret}}\mathbb{Z}[1/2].
\end{align}

\subsection{}
For any scheme $X$ we denote by $X_r$ the associated real scheme, and $D(X_r,\mathbb{Z})$ the derived category of sheaves. In the remaining of this section we study some properties of $X_r$.

\begin{lemma}
\label{lm:retmod}
If $X$ is a smooth scheme over a field, then the derived $\infty$-category $D(X_r,\mathbb{Z})$ is equivalent to the category of modules over $\mathbf{H}_{\rm{ret}}\mathbb{Z}$.
\end{lemma}
\proof
By \cite[Theorem 35]{Bac}, there is a canonical isomorphism $\mathbf{H}_{\mathbb{A}^1}\mathbb{Z}[\rho^{-1}]\simeq\mathbf{H}_{\rm{ret}}\mathbb{Z}$. Then the proof is parallel to \cite[Thm. 5.1]{CD2}, \cite[Prop. 4.3.19]{EHK}, and \cite[Cor. 5.11]{EK}.
\endproof

\subsection{}
For any scheme $X$ we denote by $C(X_r,\mathbb{Z})$ the group of continuous (or equivalently, locally constant) functions $X_r\to\mathbb{Z}$.
\begin{lemma}
\label{lm:Csurj}
Let $X$ be a scheme and let $U$ be an open subscheme of $X$. Then the canonical map $C(X_r,\mathbb{Z})\to C(U_r,\mathbb{Z})$ is surjective.
\end{lemma}
\proof
We know that any generic point of $U_r$ is still a generic point of $X_r$, and that any two points in the same connected component in $U_r$ are still in the same connected component in $X_r$. Therefore any locally constant function on $U_r$ extends to a locally constant function on $X_r$.
\endproof

\subsection{}
Recall the following result from \cite[Cor. 2.15]{Jin}: the Grothendieck group $K_0(D_c(X_r,\mathbb{Z}))$ is isomorphic to the free abelian group $Cons(X_r,\mathbb{Z})$ of constructible functions over the real scheme $X_r$ via the \emph{rank function} map
\begin{align}
\label{eq:Jin2131}
\begin{split}
K_0(D_c(X_r,\mathbb{Z}))&\xrightarrow{\sim} Cons(X_r,\mathbb{Z})\\
K&\mapsto \chi(K):x\mapsto \chi(K_x)
\end{split}
\end{align}
where $\chi(K_x)$ is the Euler characteristic of the perfect complex given by the stalk of $K$ at $x$.

\begin{lemma}
\label{lm:duallc}
A class $K\in K_0(D_c(X_r,\mathbb{Z}))$ can be represented by a dualizable object 
if and only if $\chi(K)$ is a locally constant function in $Cons(X_r,\mathbb{Z})$, that is, $\chi(K)$ is locally constant and constructible.
\end{lemma}
\proof
The proof of this claim proceeds by analogy with the \'etale case, since a sheaf is dualizable if and only if it is locally isomorphic to the constant sheaf associated to a perfect complex, see \cite[Rem. 6.3.27]{CD}.
\endproof

\begin{lemma}
\label{lm:dualext}
Let $F\in D_c(U_r,\mathbb{Z})$ be a dualizable object. Then there exists a dualizable object $G\in D_c(X_r,\mathbb{Z})$ such that $j^*G$ has the same class as $F$ in $K_0(D_c(U_r,\mathbb{Z}))$.
\end{lemma}
\proof
This is a consequence of Lemma~\ref{lm:Csurj} and Lemma~\ref{lm:duallc}.
\endproof

\section{Vanishing of generalized cohomology groups}
\label{sec:van}

\subsection{}
We use the following conventions in the rest of the paper.
Let $X$ be a connected smooth $k$-scheme. Let $i:Z\to X$ be the inclusion of a nowhere dense closed subscheme, with open complement $j:U\to X$. 

\subsection{}
\label{num:condvan}
For any integer $n$, we denote by $\mathbb{E}_n(Z/X)=[i_*\mathbbold{1}_Z[n],\mathbb{E}_X]$. Consider the following two conditions:
\begin{align}
\label{eq:H0van}
\mathbb{E}_0(Z/X)=0,
\end{align}
\begin{align}
\label{eq:H1van}
\mathbb{E}_{-1}(Z/X)=0.
\end{align}

\begin{lemma}
\label{lm:vanE0}
\begin{enumerate}
\item
If $\mathbb{E}=\mathbf{H}\mathbb{Z}$, then~\eqref{eq:H0van} and~\eqref{eq:H1van} hold.
\item
If $\mathbb{E}=\mathbf{H}_{\rm{ret}}\mathbb{Z}$ or $\mathbb{E}=\mathbf{H}_{\rm{MW}}\mathbb{Z}$, then ~\eqref{eq:H0van} holds. If further $Z$ has everywhere codimension at least $2$ in $X$, then~\eqref{eq:H1van} holds.
\end{enumerate}
\end{lemma}
In particular, when $Z$ has everywhere codimension at least $2$ in $X$, the conditions in~\ref{num:condvan} also hold for $\mathbb{E}=\mathbb{Q}_k\simeq\mathbf{H}_{\rm{MW}}\mathbb{Q}$ the rational sphere spectrum.

\proof
Indeed, one uses the \emph{niveau spectral sequence} (\cite[8.4]{DFJK})
\begin{align}
E^1_{p,q}=\oplus_{x\in Z_{(p)}}\mathbb{E}_{p+q}(\kappa(x)/X)\Rightarrow \mathbb{E}_{p+q}(Z/X).
\end{align}
Since $X$ is smooth one has $\mathbb{E}_{p+q}(\kappa(x)/X)\simeq \mathbb{E}^{-2c-p-q,-c}(\kappa(x))$, where $c$ is the codimension of $x$ in $X$. Therefore it suffices to check that $\mathbb{E}^{-2c-p-q,-c}(\kappa(x))=0$ when $p+q=0$ or $-1$. For $\mathbb{E}=\mathbf{H}\mathbb{Z}$, the vanishing follows from \cite[Lemma 3.2]{SV}; for $\mathbb{E}=\mathbf{H}_{\rm{MW}}\mathbb{Z}$ this follows from \cite[Ch. 3 Prop. 4.1.2]{BCD+}; the case $\mathbb{E}=\mathbf{H}_{\rm{ret}}\mathbb{Z}$ then follows from Scheiderer's semi-purity theorem on real \'etale cohomology (\cite[Cor. 1.12]{Sch}). 
\endproof

\begin{remark}
For $\mathbb{E}=\mathbbold{1}_k$ the integral sphere spectrum, the groups $\mathbb{E}^{-2,-1}(\kappa(x))$ and $\mathbb{E}^{-3,-2}(\kappa(x))$ are torsion but do not vanish, see \cite{RSO}.
\end{remark}

\subsection{}
We consider the following condition:
\begin{align}
\label{eq:rkext}
\textrm{The canonical map }\mathbb{E}_0(X/X)\to\mathbb{E}_0(U/U) \textrm{ is an isomorphism.}
\end{align}

\subsection{}
By the long exact sequence, if both~\eqref{eq:H0van} and~\eqref{eq:H1van} hold, then~\eqref{eq:rkext} also holds, so in particular~\eqref{eq:rkext} holds in the cases treated in Lemma~\ref{lm:vanE0}. Besides these cases, we also have the following result:
\begin{lemma}
\label{lm:H0iso}
The condition~\eqref{eq:rkext} holds when $\mathbb{E}=\mathbf{H}_{\rm{ret}}\mathbb{Z}$.
\end{lemma}
\proof
We need to show that if $X$ is an excellent regular scheme and $Z$ is a nowhere dense closed subscheme, then the following canonical map is an isomorphism:
\begin{align}
H^0(X_r,\mathbb{Z})=C(X_r,\mathbb{Z})\to H^0(U_r,\mathbb{Z})=C(U_r,\mathbb{Z}).
\end{align}
The injectivity follows from the condition~\eqref{eq:H0van}, and the surjectivity follows from Lemma~\ref{lm:Csurj}.
\endproof

\subsection{}
\label{num:RS}
Let $k$ be a field and let $\Lambda$ be a coefficient ring.
We assume that at least one of the following conditions hold:
\begin{itemize}
\item
$k$ is a perfect field which satisfies resolution of singularities;
\item
$\Lambda$ is a $\mathbb{Z}[1/p]$-algebra, where $p$ is the exponential characteristic of $k$.
\end{itemize}

\begin{lemma}
\label{lm:van2}
Assume that $Z$ is smooth over $k$, or the condition~\ref{num:RS} holds.
\begin{enumerate}
\item
If $\mathbb{E}$ is $\mathbf{H}\Lambda$ or $\mathbf{H}_{\rm{ret}}\mathbb{Z}$, then $[i_*\mathbbold{1}_Z,j_!\mathbb{E}_U]=0$ and $[\mathbbold{1}_X,j_!\mathbb{E}_U]=0$.
\item
If $\mathbb{E}=\mathbf{H}_{\rm{MW}}\Lambda$, then $[i_*\mathbbold{1}_Z,j_!\mathbb{E}_U]=0$. If further $Z$ has everywhere codimension at least $2$ in $X$, or if $2$ is invertible in $\Lambda$, then $[\mathbbold{1}_X,j_!\mathbb{E}_U]=0$.
\end{enumerate}
\end{lemma}
\proof
From the exact sequences
\begin{align}
\mathbb{E}_1(Z/Z)
\simeq
[i_*\mathbbold{1}_Z,i_*\mathbb{E}_Z[-1]]
\to
[i_*\mathbbold{1}_Z,j_!\mathbb{E}_U]
\to
\mathbb{E}_0(Z/X)=0
\end{align}
and
\begin{align}
\mathbb{E}_1(Z/Z)
\simeq
[\mathbbold{1}_X,i_*\mathbb{E}_Z[-1]]
\to
[\mathbbold{1}_X,j_!\mathbb{E}_U]
\to
\mathbb{E}_0(X/X)
\xrightarrow{\sim}
\mathbb{E}_0(U/U)
\end{align}
we see that it suffices to show that $\mathbb{E}_1(Z/Z)=0$. If $Z$ is smooth over $k$, then this follows from the fact that a sheaf does not have negative cohomology groups, see \cite[Ch. 3 Prop. 4.1.2]{BCD+} (in the case $\mathbb{E}=\mathbf{H}\mathbb{Z}$ this is a particular case of the Beilinson-Soul\'e conjecture \emph{in weight $0$}, which is well-known, see \cite{Blo}). Assuming resolution of singularities, the general case follows from the smooth case by a descent spectral sequence (see for example \cite[Thm. 3.5]{Kel}).
\endproof

\section{The localized characteristic class}
\subsection{}
\label{num:JY614}
(\cite[Remark 6.1.4]{JY})
For any morphism $g$, we denote by $g^\Delta$ the cofiber of the canonical map
\begin{align}
\label{eq:gtransform}
g^*(-)\otimes g^!\mathbbold{1}\to g^!.
\end{align}

\subsection{}
Assume that the condition in~\ref{num:RS} holds.
Let $f:X\to \operatorname{Spec}(k)$ be a smooth morphism. For any object $K\in\mathbf{SH}_c(X)$, we have the $\mathbb{E}$-valued \emph{characteristic class} of $K$ as in~\eqref{eq:CX}
\begin{align}
\label{eq:CXE}
\mathbbold{1}_X
\xrightarrow{C_X(K,\mathbb{E})}
\mathbb{E}^!_{X/k}.
\end{align}
In particular we have $C_X(\mathbbold{1}_X,\mathbb{E})=e(T_{X/k},\mathbb{E})$ is the $\mathbb{E}$-valued Euler class of the tangent bundle of $X$.
By \cite[\S 5.2]{AS}, if we denote $\delta:X\to X\times X$ the diagonal morphism, then the cofiber of the map~\eqref{eq:CXE} is isomorphic to $\delta^\Delta\delta_*\mathbb{E}^!_{X/k}$, so that we have a cofiber sequence
\begin{align}
\label{eq:eucof}
\mathbbold{1}_X
\xrightarrow{C_X(K,\mathbb{E})}
\mathbb{E}^!_{X/k}
\xrightarrow{}
\delta^\Delta\delta_*\mathbb{E}^!_{X/k}.
\end{align}

\begin{remark}
From the Euler fiber sequence (\cite[Def. 3.1.6]{DJK})
\begin{align}
\Sigma^\infty_X E^\times\to \mathbbold{1}_X\xrightarrow{e(E)} Th(E)
\end{align}
for every vector bundle $E$ over $X$, we see that there exists a canonical isomorphism $\Sigma^\infty_XT_{X/k}^\times[1]\simeq\delta^\Delta\delta_*\mathcal{K}_{X/k}$.
\end{remark}

\begin{lemma}[\textrm{\cite[Prop. 4.2.6]{JY}}]
\label{lm:JY426}
If $L\to M\to N$ is a distinguished triangle in $\mathbf{SH}_c(X)$, then one has $C_X(M,\mathbb{E})=C_X(L,\mathbb{E})+C_X(N,\mathbb{E})$.
\end{lemma}

\subsection{}
In what follows, we let $K\in\mathbf{SH}_c(X)$ be such that $K_{|U}$ is dualizable. 

\begin{definition}
\label{def:rankE}
\begin{enumerate}
\item
We denote by $\chi(K_{|U},\mathbb{E})$ the $\mathbb{E}$-valued Euler characteristic of $K_{|U}$, that is, the trace of its identity map, considered as an element of $\mathbb{E}_0(U/U)$. 
\item
If the condition~\eqref{eq:rkext} holds, we define the ($\mathbb{E}$-valued) \textbf{rank} of $K$, denoted as $\operatorname{rk}(K,\mathbb{E})$, as the unique element in $\mathbb{E}_0(X/X)$ whose image in $\mathbb{E}_0(U/U)$ is $\operatorname{rk}(K_{|U},\mathbb{E})$.
\end{enumerate}
\end{definition}

\subsection{}
\label{num:rkwd}
By Lemma~\ref{lm:vanE0} and Lemma~\ref{lm:H0iso}, the rank is well-defined in the following cases: 
\begin{enumerate}
\item
$\mathbb{E}=\mathbf{H}\mathbb{Z}$ or $\mathbf{H}_{\rm{ret}}\mathbb{Z}$,
\item
$\mathbb{E}=\mathbf{H}_{\rm{MW}}\Lambda$ or $\mathbf{H}_{\rm{W}}\mathbb{Z}$, with either $2$ invertible in $\Lambda$, or $Z$ everywhere of codimension at least $2$ in $X$.
\end{enumerate}

\begin{lemma}
\label{lm:addrk}
In the setting of~\ref{num:rkwd}, if $L\to M\to N$ is a distinguished triangle of dualizable objects in $\mathbf{SH}_c(X)$, then one has $\operatorname{rk}(M,\mathbb{E})=\operatorname{rk}(L,\mathbb{E})+\operatorname{rk}(N,\mathbb{E})$.
\end{lemma}
\proof
This follows from the additivity of Euler characteristics, see \cite{May}.
\endproof

\begin{lemma}
The following diagram is commutative:
\begin{align}
\begin{split}
  \xymatrix{
    j_!\mathbbold{1}_U\ar[r]^-{\chi(K_{|U},\mathbb{E})} \ar[d]_-{} & \mathbb{E}_X \ar[d]^-{e(T_{X/k},\mathbb{E})} \\
    \mathbbold{1}_X\ar[r]_-{C_X(K,\mathbb{E})} & \mathbb{E}^!_{X/k}
    }
  \end{split}
\end{align}
\end{lemma}
\proof
See \cite[Lemma 5.1.10]{JY}.
\endproof

\subsection{}
We define the map
\begin{align}
\label{eq:mapZ}
B_X^Z(K,\mathbb{E}):
i_*\mathbbold{1}_Z
\to
\delta^\Delta\delta_* \mathbb{E}^!_{X/k}
\end{align}
as the composition
\begin{align}
i_*\mathbbold{1}_Z
\to
i_*i^*\delta^\Delta\delta_*\mathbbold{1}_X
\to
i_*i^*\delta^\Delta(\mathbb{D}(K)\boxtimes_kK)
\simeq
\delta^\Delta(\mathbb{D}(K)\boxtimes_kK)
\to
\delta^\Delta\delta_* \mathbb{E}^!_{X/k}.
\end{align}

\begin{lemma}
\label{lm:dualb0}
If $K$ is a dualizable object, then $B_X^Z(K,\mathbb{E})=0$.
\end{lemma}
\proof
This is because the transformation~\eqref{eq:gtransform} is an isomorphism when applied to dualizable objects (\cite[5.4]{FHM}), and therefore $\delta^\Delta(\mathbb{D}(K)\boxtimes_kK)=0$.
\endproof

\begin{lemma}
\label{lm:addB}
If $L\to M\to N$ is a distinguished triangle in $\mathbf{SH}_c(X)$, then one has $B_X^Z(M,\mathbb{E})=B_X^Z(L,\mathbb{E})+B_X^Z(N,\mathbb{E})$.
\end{lemma}

\proof
The claim is a form of additivity of Euler characteristics, which can be proved in the style of \cite{May}, \cite{GPS}. We proceed as in the proof of \cite[Prop. 4.2.6]{JY}: using the language of higher categories, there exist objects $u,v,w$ and a commutative diagram of the form
\begin{align}
\begin{split}
\xymatrix{
    & i_*i^*\delta^\Delta\delta_*\mathbbold{1}_X \ar^-{}[d] \ar^-{}[rd] \ar_-{}[ld] &\\
   i_*i^*\delta^\Delta(\mathbb{D}(L)\boxtimes_kL\oplus\mathbb{D}(N)\boxtimes_kN) \ar^-{\wr}[dd] \ar^-{}[rd]
   & i_*i^*\delta^\Delta\delta_*v \ar^-{}[r] \ar^-{}[l] 
   & i_*i^*\delta^\Delta(\mathbb{D}(M)\boxtimes_kM) \ar^-{}[ld] \ar_-{\wr}[dd]\\
   & i_*i^*\delta^\Delta\delta_*w \ar@{.>}^-{}[d]
   &\\
   \delta^\Delta(\mathbb{D}(L)\boxtimes_kL\oplus\mathbb{D}(N)\boxtimes_kN) \ar^-{}[r]  \ar^-{}[rd]
   & \delta^\Delta u \ar^-{}[d]
   & \delta^\Delta(\mathbb{D}(M)\boxtimes_kM) \ar^-{}[l] \ar^-{}[ld]\\
   & \delta^\Delta\delta_* \mathbb{E}^!_{X/k} &
 }
 \end{split}
\end{align}
which finishes the proof.
\endproof

\begin{definition}
\label{def:swan}
Let $K\in\mathbf{SH}_c(X)$ be such that $K_{|U}$ is dualizable.
Assume that the condition~\ref{num:RS} holds. The \textbf{($\mathbb{E}$-valued) localized characteristic class}
, if it exists, is the unique class 
\begin{align}
C_X^Z(K,\mathbb{E})\in\mathbb{E}_0(Z/k)\simeq[i_*\mathbbold{1}_Z, \mathbb{E}^!_{X/k}]
\end{align}
whose image under the canonical map 
\begin{align}
\label{eq:phican}
\phi_X^Z:
[i_*\mathbbold{1}_Z, \mathbb{E}^!_{X/k}]
\to
[i_*\mathbbold{1}_Z, \delta^\Delta\delta_*\mathbb{E}^!_{X/k}]
\end{align}
is the class 
$B_X^Z(K,\mathbb{E})\in [i_*\mathbbold{1}_Z,
\delta^\Delta\delta_*\mathbb{E}^!_{X/k}]$ in~\eqref{eq:mapZ}. Here the uniqueness is a consequence of the vanishing~\eqref{eq:H0van}.

\end{definition}

\begin{lemma}
\label{lm:locprop}
\begin{enumerate}
\item
If $K$ is a dualizable object, then $C_X^Z(K,\mathbb{E})=0$.
\item
If $L\to M\to N$ is a distinguished triangle in $\mathbf{SH}_c(X)$, then one has $C_X^Z(M,\mathbb{E})=C_X^Z(L,\mathbb{E})+C_X^Z(N,\mathbb{E})$, if all these classes exist.
\end{enumerate}
\end{lemma}
\proof
\begin{enumerate}
\item
Follows from Lemma~\ref{lm:dualb0}.
\item
Follows from Lemma~\ref{lm:addB}.
\end{enumerate}
\endproof

\begin{lemma}
\label{lm:iCZ}
If $K=i_*F$, then the class $C_X^Z(K,\mathbb{E})$ exists and satisfies the formula
\begin{align}
C_X^Z(i_*F,\mathbb{E})=C_Z(F,\mathbb{E})\in\mathbb{E}_0(Z/k).
\end{align}
\end{lemma}
\proof
We need to show that $\phi_X^Z(C_Z(F,\mathbb{E}))=B_X^Z(i_*F,\mathbb{E})$, where $\phi_X^Z$ is the map in~\eqref{eq:phican}.
Using the isomorphism $\mathbb{D}(i_*F)\boxtimes_ki_*F\simeq (i\times i)_*(\mathbb{D}(F)\boxtimes_kF)$, the result follows from the following commutative diagram:
\begin{align}
\begin{split}
  \xymatrix{
     i_*\mathbbold{1}_Z \ar[r]^-{} \ar[d]^-{} & i_*(\mathbb{D}(F)\otimes_kF) \ar[r]^-{} \ar[d]^-{} & i_*\mathbb{E}^!_{Z/k} \ar[r]^-{}  \ar[d]_-{} & \mathbb{E}^!_{X/k} \ar[d]^-{} \\ 
    i_*i^*\delta^\Delta\delta_*\mathbbold{1}_X \ar[d]^-{} & i_*\delta_Z^!(\mathbb{D}(F)\boxtimes_kF) \ar[d]_-{} \ar[ru]^-{} & i_*\delta_Z^\Delta\delta_{Z*}\mathbb{E}^!_{Z/k}\ar[r]_-{} & \delta^\Delta\delta_*\mathbb{E}^!_{X/k} \\
    i_*i^*\delta^\Delta(\mathbb{D}(i_*F)\boxtimes_ki_*F) \ar[r]^-{\sim} & i_*\delta_Z^\Delta(\mathbb{D}(F)\boxtimes_kF) \ar[r]^-{\sim} \ar[ru]^-{} & \delta^\Delta(\mathbb{D}(i_*F)\boxtimes_ki_*F). \ar[ru]^-{} &  
  }
  \end{split}
\end{align}
\endproof

\begin{proposition}
\label{prop:locex}
The class $C_X^Z(K,\mathbb{E})$ exists in the following cases:
\begin{enumerate}
\item
If there exists a dualizable object $G$ such that $j^*K\simeq j^*G$;
\item 
The spectrum $\mathbb{E}$ is either $\mathbf{H}\Lambda$ or $\mathbf{H}_{\rm{ret}}\mathbb{Z}$;
\item
The spectrum $\mathbb{E}$ is either $\mathbb{E}=\mathbf{H}_{\rm{MW}}\Lambda$ or $\mathbb{E}=\mathbf{H}_{\rm{W}}\Lambda$, and at least one of the following conditions is satisfied:
\begin{enumerate}
\item
$2$ invertible in $\Lambda$;
\item
$Z$ everywhere of codimension at least $2$ in $X$;
\item
$X$ is odd-dimensional.
\end{enumerate}
\end{enumerate}
\end{proposition}
\proof
We first consider the case where there exists a dualizable object $G$ such that $j^*K\simeq j^*G$. In this case, we define the class $C_X^Z(K,\mathbb{E})$ by the formula
\begin{align}
C_X^Z(K,\mathbb{E})
=
C_Z(i^*K,\mathbb{E})
-
C_Z(i^*G,\mathbb{E}).
\end{align}
By additivity (Lemma~\ref{lm:addB}) and Lemma~\ref{lm:iCZ}, the result follows from the following equalities:
\begin{align}
\phi_X^Z(C_Z(i^*K,\mathbb{E}))=B_X^Z(i_*i^*K,\mathbb{E})
\end{align}
\begin{align}
\phi_X^Z(-C_Z(i^*G,\mathbb{E}))=-B_X^Z(i_*i^*G,\mathbb{E})=B_X^Z(j_!j^*G,\mathbb{E})=B_X^Z(j_!j^*K,\mathbb{E})
\end{align}
where we use the fact that $B_X^Z(G,\mathbb{E})=0$ by Lemma~\ref{lm:dualb0}.

For the remaining cases, if $\mathbb{E}=\mathbf{H}\Lambda$, or if $\mathbb{E}=\mathbf{H}_{\rm{MW}}\Lambda$ and $Z$ has everywhere of codimension at least $2$, then the existence follows from the vanishing~\eqref{eq:H1van}. The case $\mathbb{E}=\mathbf{H}_{\rm{MW}}\Lambda$ with $2$ invertible in $\Lambda$ would follow from the cases $\mathbb{E}=\mathbf{H}\Lambda$ and $\mathbb{E}=\mathbf{H}_{\rm{ret}}\mathbb{Z}$ using the isomorphism~\eqref{eq:MW2split}. 
Therefore it remains to show the existence in the case $\mathbb{E}=\mathbf{H}_{\rm{ret}}\mathbb{Z}$, or $\mathbb{E}=\mathbf{H}_{\rm{MW}}\Lambda$ and $X$ odd-dimensional.

Consider the case $\mathbb{E}=\mathbf{H}_{\rm{ret}}\mathbb{Z}$. Since $\mathbf{H}_{\rm{ret}}\mathbb{Z}$ is a commutative ring spectrum with unit, the commutative diagram
\begin{align}
\begin{split}
  \xymatrix{
   i_*i^*\delta^\Delta\delta_*\mathbbold{1}_X \ar[r]^-{} \ar[d]^-{} & i_*i^*\delta^\Delta(\mathbb{D}(K\otimes\mathbf{H}_{\rm{ret}}\mathbb{Z})\boxtimes_k(K\otimes\mathbf{H}_{\rm{ret}}\mathbb{Z}))  \ar[d]^-{\wr}  \\ 
    i_*i^*\delta^\Delta(\mathbb{D}(K)\boxtimes_kK) \ar[d]^-{\wr} \ar[ru]^-{} & \delta^\Delta(\mathbb{D}(K\otimes\mathbf{H}_{\rm{ret}}\mathbb{Z})\boxtimes_k(K\otimes\mathbf{H}_{\rm{ret}}\mathbb{Z})) \ar[d]^-{} \\   
  \delta^\Delta(\mathbb{D}(K)\boxtimes_kK) \ar[ru]^-{} \ar[r]^-{} & \delta^\Delta\delta_*\mathbf{H}_{\rm{ret}}\mathbb{Z}^!_{X/k} 
  }
  \end{split}
\end{align}
gives the identity
\begin{align}
B_X^Z(K,\mathbf{H}_{\rm{ret}}\mathbb{Z})
=
B_X^Z(K\otimes\mathbf{H}_{\rm{ret}}\mathbb{Z},\mathbf{H}_{\rm{ret}}\mathbb{Z}).
\end{align}
The object $K\otimes\mathbf{H}_{\rm{ret}}\mathbb{Z}$ lives in $D_c(X_r,\mathbb{Z})$ by Lemma~\ref{lm:retmod}, and since the functor
\begin{align}
\begin{split}
\mathbf{SH}&\to D(X_r,\mathbb{Z})\\
K&\mapsto K\otimes\mathbf{H}_{\rm{ret}}\mathbb{Z}
\end{split}
\end{align}
together with the forgetful functor $D(X_r,\mathbb{Z})\to\mathbf{SH}$ form a premotivic adjunction, we may perform the construction of the class $B_X^Z(K\otimes\mathbf{H}_{\rm{ret}}\mathbb{Z},\mathbf{H}_{\rm{ret}}\mathbb{Z})$ in the motivic category $D_c(X_r,\mathbb{Z})$, which still yields the same class.

In the category $D_c(X_r,\mathbb{Z})$, by Lemma~\ref{lm:dualext}, there exists a dualizable object $G\in D_c(X_r,\mathbb{Z})$ such that $j^*G$ has the same class as $j^*K$ in $K_0(D_c(U_r,\mathbb{Z}))$. Then by additivity, as the case above, the class $C_X^Z(K,\mathbb{E})$ exists and satisfies the formula
\begin{align}
C_X^Z(K,\mathbf{H}_{\rm{ret}}\mathbb{Z})
=
C_Z(i^*K,\mathbf{H}_{\rm{ret}}\mathbb{Z})
-
C_Z(i^*G,\mathbf{H}_{\rm{ret}}\mathbb{Z})
\end{align}
which is well-defined and does not depend on the choice of $G$. This proves the existence in the case $\mathbb{E}=\mathbf{H}_{\rm{ret}}\mathbb{Z}$.

Now assume that $\mathbb{E}=\mathbf{H}_{\rm{MW}}\Lambda$ or $\mathbb{E}=\mathbf{H}_{\rm{W}}\Lambda$ and $X$ is odd-dimensional. We know that the group $[i_*\mathbbold{1}_Z, \delta^\Delta\delta_*\mathbb{E}^!_{X/k}]$ fits into a short exact sequence
\begin{align}
\label{eq:deltases}
0\to \mathbb{E}_0(Z/k)\to
[i_*\mathbbold{1}_Z, \delta^\Delta\delta_*\mathbb{E}^!_{X/k}]
\to \mathbb{E}_{-1}(Z/X)\to0.
\end{align}
and therefore the element $B_X^Z(K,\mathbb{E})$ lifts to $\mathbb{E}_0(Z/k)$ if and only if its image in $\mathbb{E}_{-1}(Z/X)$ is zero. We know that the canonical map $\mathbf{H}_{\rm{MW}}\Lambda\to\mathbf{H}_{\rm{W}}\Lambda$ induces an isomorphism 
\begin{align}
\mathbf{H}_{\rm{MW}}\Lambda_{-1}(Z/X)\simeq\mathbf{H}_{\rm{W}}\Lambda_{-1}(Z/X),
\end{align}
and therefore it suffices to show that for $\mathbb{E}=\mathbf{H}_{\rm{W}}\Lambda$ the element $B_X^Z(K,\mathbf{H}_{\rm{W}}\Lambda)$ has a lift. If $X$ is odd-dimensional, then its tangent bundle has odd rank, and therefore the Witt-valued Euler class $e(T_{X/k},\mathbf{H}_{\rm{W}}\Lambda)$ is trivial (\cite[Lemma 4.3]{LevWitt}, \cite[Prop. 4.4]{LevA}). Therefore the cofiber sequence~\eqref{eq:eucof} splits, and the canonical map $\mathbf{H}_{\rm{W}}\Lambda^!_{X/k}\to\delta^\Delta\delta_*\mathbf{H}_{\rm{W}}\Lambda^!_{X/k}$ has a retract 
\begin{align}
\label{eq:Eulerret}
\delta^\Delta\delta_*\mathbf{H}_{\rm{W}}\Lambda^!_{X/k}\to\mathbf{H}_{\rm{W}}\Lambda^!_{X/k}
\end{align}
and we obtain a lift of the element $B_X^Z(K,\mathbf{H}_{\rm{W}}\Lambda)$ to $\mathbf{H}_{\rm{W}}\Lambda_0(Z/k)$ by composing it with the map~\eqref{eq:Eulerret}, which finishes the proof.
\endproof



\section{A Grothendieck-Ogg-Shafarevich-type formula}

\subsection{}
Let $K\in\mathbf{SH}_c(X)$ be such that $K_{|U}$ is dualizable. In the situation of~\ref{num:rkwd}, the $\mathbb{E}$-valued rank $\operatorname{rk}(K,\mathbb{E})$ as in Definition~\ref{def:rankE} is a well-defined element in $\mathbb{E}_0(X/X)$. 

\begin{theorem}
\label{th:GOS}
Assume that the condition~\ref{num:RS} holds, and that in addition $\mathbb{E}$ is one of the following: 
\begin{enumerate}
\item
$\mathbb{E}=\mathbf{H}\Lambda$; 
\item
$\mathbb{E}=\mathbf{H}_{\rm{ret}}\mathbb{Z}$; 
\item
$\mathbb{E}=\mathbf{H}_{\rm{MW}}\Lambda$ or $\mathbf{H}_{\rm{W}}\Lambda$, with
\begin{enumerate}
\item either $2$ is invertible in $\Lambda$,
\item or $Z$ everywhere of codimension at least $2$ in $X$, and 
$X$ is odd-dimensional.
\end{enumerate}
\end{enumerate}
Then one has
\begin{align}
\label{eq:GOS}
C_X(K,\mathbb{E})=\operatorname{rk}(K,\mathbb{E})\cdot C_X(\mathbbold{1}_X,\mathbb{E})+i_*C_X^Z(K,\mathbb{E})\in \mathbb{E}_0(X/k).
\end{align}
\end{theorem}

\begin{remark}
In the case $\mathbb{E}=\mathbf{H}_{\rm{MW}}\Lambda$ or $\mathbf{H}_{\rm{W}}\Lambda$, we believe that the result also holds when $Z$ has everywhere codimension at least $2$ in $X$ and $X$ even-dimensional, but we do not know how to prove it.
\end{remark}

\subsection{}
\label{num:Wittcase}
We begin the proof of Theorem~\ref{th:GOS}.
First assume that $\mathbb{E}=\mathbf{H}_{\rm{W}}\Lambda$ and $X$ is odd-dimensional. As in the proof of Proposition~\ref{prop:locex}, since the Witt-valued Euler class $e(T_{X/k},\mathbf{H}_{\rm{W}}\Lambda)$ is trivial (\cite[Lemma 4.3]{LevWitt}, \cite[Prop. 4.4]{LevA}), we have a commutative diagram
\begin{align}
\begin{split}
  \xymatrix{
    \mathbbold{1}_X \ar[rr]^-{C_X(K,\mathbf{H}_{\rm{W}}\Lambda)} \ar[d]_-{} & & \mathbf{H}_{\rm{W}}\Lambda^!_{X/k} \\
    i_*\mathbbold{1}_Z \ar[rr]_-{\eqref{eq:mapZ}} \ar[rru]^-{C_X^Z(K,\mathbf{H}_{\rm{W}}\Lambda)} & & \delta^\Delta\delta_*\mathbf{H}_{\rm{W}}\Lambda^!_{X/k} \ar[u]_-{\eqref{eq:Eulerret}}
    }
  \end{split}
\end{align}
which shows that
\begin{align}
\label{eq:CZXW}
C_X(K,\mathbf{H}_{\rm{W}}\Lambda)=i_*C_X^Z(K,\mathbf{H}_{\rm{W}}\Lambda)\in \mathbf{H}_{\rm{W}}\Lambda_0(X/k).
\end{align}
Since $C_X(\mathbbold{1}_X,\mathbb{E})=e(T_{X/k},\mathbf{H}_{\rm{W}}\Lambda)=0$, this proves the formula~\eqref{eq:GOS}.

\subsection{}
By Lemma~\ref{lm:JY426}, Lemma~\ref{lm:addrk} and Lemma~\ref{lm:locprop}, all the three terms of the formula~\eqref{eq:GOS} are additive along distinguished triangles. Therefore by localization it suffices to prove the following two formulas:
\begin{align}
\label{eq:GOSiF}
C_X(i_*i^*K,\mathbb{E})=\operatorname{rk}(i_*i^*K,\mathbb{E})\cdot C_X(\mathbbold{1}_X,\mathbb{E})+i_*C_X^Z(i_*i^*K,\mathbb{E})\in \mathbb{E}_0(X/k).
\end{align}
\begin{align}
\label{eq:GOSjF}
C_X(j_!j^*K,\mathbb{E})=\operatorname{rk}(j_!j^*K,\mathbb{E})\cdot C_X(\mathbbold{1}_X,\mathbb{E})+i_*C_X^Z(j_!j^*K,\mathbb{E})\in \mathbb{E}_0(X/k).
\end{align}
For the formula~\eqref{eq:GOSiF}, since $\operatorname{rk}(i_*i^*K,\mathbb{E})=0$, the result follows from Lemma~\ref{lm:iCZ}. Therefore it remains to prove the formula~\eqref{eq:GOSjF}, which amounts to say that for every dualizable object $F\in\mathbf{SH}_c(U)$ one has
\begin{align}
\label{eq:GOSF}
C_X(j_!F,\mathbb{E})=\operatorname{rk}(j_!F,\mathbb{E})\cdot C_X(\mathbbold{1}_X,\mathbb{E})+i_*C_X^Z(j_!F,\mathbb{E})\in \mathbb{E}_0(X/k).
\end{align}

\subsection{}
In the case $\mathbb{E}=\mathbf{H}_{\rm{ret}}\mathbb{Z}$, we may proceed as in the proof of Proposition~\ref{prop:locex}: we may replace $F$ by $F\otimes\mathbf{H}_{\rm{ret}}\mathbb{Z}$ and therefore work in the category $D_c(X_r,\mathbb{Z})$. By Lemma~\ref{lm:dualext}, there exists a dualizable object $G\in D_c(X_r,\mathbb{Z})$ such that $j^*G$ has the same class as $F$ in $K_0(D_c(U_r,\mathbb{Z}))$. Using Lemma~\ref{lm:locprop} and Lemma~\ref{lm:iCZ}, we obtain
\begin{align}
\begin{split}
C_X(j_!F,\mathbf{H}_{\rm{ret}}\mathbb{Z})
&=
C_X(j_!j^*G,\mathbf{H}_{\rm{ret}}\mathbb{Z})
=
C_X(G,\mathbf{H}_{\rm{ret}}\mathbb{Z})-i_*C_Z(i^*G,\mathbf{H}_{\rm{ret}}\mathbb{Z})\\
&=
\operatorname{rk}(j_!F,\mathbf{H}_{\rm{ret}}\mathbb{Z})\cdot C_X(\mathbbold{1}_X,\mathbf{H}_{\rm{ret}}\mathbb{Z})-i_*C_X^Z(i_*i^*G,\mathbf{H}_{\rm{ret}}\mathbb{Z})\\
&=
\operatorname{rk}(j_!F,\mathbf{H}_{\rm{ret}}\mathbb{Z})\cdot C_X(\mathbbold{1}_X,\mathbf{H}_{\rm{ret}}\mathbb{Z})-i_*C_X^Z(G,\mathbf{H}_{\rm{ret}}\mathbb{Z})+i_*C_X^Z(j_!j^*G,\mathbf{H}_{\rm{ret}}\mathbb{Z})\\
&=
\operatorname{rk}(j_!F,\mathbf{H}_{\rm{ret}}\mathbb{Z})\cdot C_X(\mathbbold{1}_X,\mathbf{H}_{\rm{ret}}\mathbb{Z})+i_*C_X^Z(j_!F,\mathbf{H}_{\rm{ret}}\mathbb{Z}).
\end{split}
\end{align}
which proves~\eqref{eq:GOSF}.

\subsection{}
We now prove the remaining cases.
Without loss of generality, we may assume that $X$, and therefore $U$, is connected.
In general, it suffices to show the following equality:
\begin{align}
\label{eq:CZtr}
C_X(j_!F,\mathbb{E})-\operatorname{rk}(F,\mathbb{E})\cdot C_X(j_!\mathbbold{1}_U,\mathbb{E})
=
i_*C_X^Z(j_!F,\mathbb{E})-\operatorname{rk}(F,\mathbb{E})\cdot i_*C_X^Z(j_!\mathbbold{1}_U,\mathbb{E}).
\end{align}

\subsection{}
We define the following two maps:
\begin{align}
\begin{split}
C_{X,U}(F,\mathbb{E}):
\mathbbold{1}_X
&\to
\underline{Hom}(j_!F,j_!F)
\simeq
\delta^!(\mathbb{D}(j_!F)\boxtimes_kj_!F)\\
&\to
\delta^*(\mathbb{D}(j_!F)\boxtimes_kj_!F)
\simeq
j_!(\mathbb{D}(F)\otimes_kF)
\to
j_!\mathbb{E}^!_{U/k},
\end{split}
\end{align}
\begin{align}
\begin{split}
C_{X,U}^Z(F,\mathbb{E})&:
i_*\mathbbold{1}_Z
\to
i_*i^*\delta^\Delta\delta_*\mathbbold{1}_X
\to
i_*i^*\delta^\Delta\delta_*\underline{Hom}(j_!F,j_!F)\\
&\simeq
i_*i^*\delta^\Delta\delta_*\delta^!(\mathbb{D}(j_!F)\boxtimes_kj_!F)
\to
i_*i^*\delta^\Delta(\mathbb{D}(j_!F)\boxtimes_kj_!F)\\
&\simeq
\delta^\Delta(\mathbb{D}(j_!F)\boxtimes_kj_!F)
\to
\delta^\Delta\delta_*j_!(\mathbb{D}(F)\otimes_kF)
\to
\delta^\Delta\delta_*j_!\mathbb{E}^!_{U/k}.
\end{split}
\end{align}

\begin{lemma}
\label{lm:openCX}
\begin{enumerate}
\item
The composition
\begin{align}
\mathbbold{1}_X\xrightarrow{C_{X,U}(F,\mathbb{E})}j_!\mathbb{E}^!_{U/k}\to\mathbb{E}^!_{X/k}
\end{align} 
agrees with $C_X(j_!F,\mathbb{E})$.
\item
The composition 
\begin{align}
i_*\mathbbold{1}_Z
\xrightarrow{C_{X,U}^Z(F,\mathbb{E})}
\delta^\Delta\delta_*j_!\mathbb{E}^!_{U/k}
\to
\delta^\Delta\delta_*\mathbb{E}^!_{X/k}
\end{align}
agrees with the class $B_X^Z(j_!F,\mathbb{E})$ in~\eqref{eq:mapZ}.
\end{enumerate}
\end{lemma}
\proof
Both claims follow from the following commutative diagram:
\begin{align}
\begin{split}
  \xymatrix@=10pt{
    \mathbb{D}(j_!F)\otimes_kj_!F \ar[r]^-{\sim} \ar[d]_-{} & j_!(\mathbb{D}(F)\otimes_kF) \ar[d]_-{} \\
    \mathbb{E}^!_{X/k} & j_!\mathbb{E}^!_{U/k}. \ar[l]_-{}
    }
  \end{split}
\end{align}
\endproof

\subsection{}
We have the following commutative diagram:
\begin{align}
\label{eq:ASdiag}
\begin{split}
  \xymatrix{
      & [i_*\mathbbold{1}_Z,j_!\mathbb{E}_U]=0 \ar@{=}[r] \ar[d]^-{} & [\mathbbold{1}_X,j_!\mathbb{E}_U] \ar[d]_-{} \\ 
      & [i_*\mathbbold{1}_Z,j_!\mathbb{E}^!_{U/k}] \ar[d]_-{a} \ar[r]^-{c} & [\mathbbold{1}_X,j_!\mathbb{E}^!_{U/k}] \ar[d]_-{b} \\
    [\mathbbold{1}_X,\mathbb{E}_X] \ar[d]_-{\wr} & [i_*\mathbbold{1}_Z,\delta^\Delta\delta_*j_!\mathbb{E}^!_{U/k}] \ar[r]^-{} \ar[d]_-{\partial} & [\mathbbold{1}_X,\delta^\Delta\delta_*j_!\mathbb{E}^!_{U/k}] \ar[d]_-{} \\ 
  [\mathbbold{1}_U,\mathbb{E}_U] \ar[r]^-{d} & [i_*\mathbbold{1}_Z,j_!\mathbb{E}_U[1]] \ar[r]^-{} & [\mathbbold{1}_X,j_!\mathbb{E}_U[1]].
  }
  \end{split}
\end{align}
Since both columns are exact, we know that both maps $a$ and $b$ are injective.
Consider the following two classes:
\begin{align}
\alpha=
C_{X,U}(F,\mathbb{E})-\operatorname{rk}(F,\mathbb{E})\cdot C_{X,U}(\mathbbold{1}_U,\mathbb{E})
\in 
[\mathbbold{1}_X,j_!\mathbb{E}^!_{U/k}]
\end{align}
and
\begin{align}
\beta=
C_{X,U}^Z(F,\mathbb{E})-\operatorname{rk}(F,\mathbb{E})\cdot C_{X,U}^Z(\mathbbold{1}_U,\mathbb{E})
\in 
[i_*\mathbbold{1}_Z,\delta^\Delta\delta_*j_!\mathbb{E}^!_{U/k}].
\end{align}
By Lemma~\ref{lm:openCX}, we know that 
\begin{enumerate}
\item
the image of $\alpha$ by the map $[\mathbbold{1}_X,j_!\mathbb{E}^!_{U/k}]\to [\mathbbold{1}_X,\mathbb{E}^!_{X/k}]$ is $C_X(j_!F,\mathbb{E})-\operatorname{rk}(F,\mathbb{E})\cdot C_X(j_!\mathbbold{1}_U,\mathbb{E})$,
\item
the image of $\beta$ by the isomorphism $[i_*\mathbbold{1}_Z,\delta^\Delta\delta_*\mathbb{E}^!_{X/k}]\simeq[i_*\mathbbold{1}_Z,\mathbb{E}^!_{X/k}]$ is $C_X^Z(j_!F,\mathbb{E})-\operatorname{rk}(F,\mathbb{E})\cdot C_X^Z(j_!\mathbbold{1}_U,\mathbb{E})$. 
\end{enumerate}
Therefore in order to show the equality~\eqref{eq:CZtr}, we only need to show that there exists a class $\gamma\in[i_*\mathbbold{1}_Z,j_!\mathbb{E}^!_{U/k}]$ such that $a(\gamma)=\beta$ and $c(\gamma)=\alpha$. Since the two classes $\alpha$ and $\beta$ have the same image in $[\mathbbold{1}_X,\delta^\Delta\delta_*j_!\mathbb{E}^!_{U/k}]$, by the commutative diagram~\eqref{eq:ASdiag} we only need to show that $\partial(\beta)=0\in[i_*\mathbbold{1}_Z,j_!\mathbb{E}_U[1]]$. Using the diagram~\eqref{eq:ASdiag} again, we know that there exists a class $\sigma\in[\mathbbold{1}_U,\mathbb{E}_U]\simeq[\mathbbold{1}_X,\mathbb{E}_X]$ such that $d(\sigma)=\partial(\beta)$. By Lemma~\ref{lm:van2}, we know that the map $d$ is injective, and therefore we are reduced to show that
\begin{align}
\label{eq:sigma0}
\sigma=0\in[\mathbbold{1}_X,\mathbb{E}_X]=\mathbb{E}_0(X/X).
\end{align}

\subsection{}
In the case $\mathbb{E}=\mathbf{H}\Lambda$, we have $\sigma\in[\mathbbold{1}_X,\mathbbold{1}_X]=\Lambda$. Since the six functors, and therefore the construction of the class $\sigma$ is compatible with \'etale realizations, we know that the image of $\sigma$ in $\mathbb{Z}/\ell\mathbb{Z}$ is $0$ for every prime $\ell$ different from $p$ (see the proof of \cite[Lemma 5.2.4]{AS}). If $\sigma\ne0\in\Lambda$, one may take a sufficiently large prime $\ell$ to obtain a contradiction, which proves the claim.

\subsection{}
In the case $\mathbb{E}=\mathbf{H}_{\rm{MW}}\Lambda$, since the canonical map
\begin{align}
\mathbf{H}_{\rm{MW}}\Lambda_0(X/X)\to \mathbf{H}_{\rm{W}}\Lambda_0(X/X)\times \mathbf{H}_{\rm{MW}}\Lambda_0(X/X)
\end{align}
is injective, it remains to check the equality~\eqref{eq:sigma0} for $\mathbb{E}=\mathbf{H}_{\rm{W}}\Lambda$. If $X$ is odd-dimensional, this is proven in~\ref{num:Wittcase}. If $2$ is invertible in $\Lambda$, the result follows from the isomorphism~\eqref{eq:MW2split} and the case $\mathbb{E}=\mathbf{H}_{\rm{ret}}\mathbb{Z}$ proved above. Therefore we have finished the proof of Theorem~\ref{th:GOS}.

\begin{corollary}
\label{cor:GOS2p}
In the setting of Theorem~\ref{th:GOS}, further assume that the structure morphism $p:X\to \operatorname{Spec}(k)$ is proper (and smooth). Then one has
\begin{align}
\label{eq:pfgos}
\chi(p_*K)=p_*(\operatorname{rk}(K,\mathbf{H}_{\rm{MW}}\Lambda)\cdot e(T_{X/k},\mathbf{H}_{\rm{MW}}\Lambda))+\int_ZC_X^Z(K,\mathbf{H}_{\rm{MW}}\Lambda)\in GW(k)\otimes\Lambda.
\end{align}
\end{corollary}

\subsection{}
If the exponential characteristic $p$ of $k$ is different from $2$, by \cite[Rem. 2.1]{LevA}, there is a Cartesian diagram of the form
\begin{align}
\label{eq:GWp}
\begin{split}
  \xymatrix@=10pt{
    GW(k) \ar[r]^-{} \ar[d]_-{} & \mathbb{Z} \ar[d]_-{} \\
    GW(k)[1/p] \ar[r]_-{} & \mathbb{Z}[1/p]. 
    }
  \end{split}
\end{align}

\begin{definition}
\label{def:quadSwan}
Let $Art(K_{\textrm{et}})\in\mathbb{Z}$ be the Artin conductor defined using the \'etale realization of $K$ as
\begin{align}
Art(K_{\textrm{et}})=\int Sw(K_{\textrm{et}|U})+\operatorname{rk}(K_{\textrm{et}})\cdot \chi(i^*\Lambda)-\chi(i^*K_{\textrm{et}})\in\mathbb{Z}
\end{align}
where $\int Sw(K_{\textrm{et}|U})$ is the Swan conductor of the \'etale realization of $K_{|U}$, as appears in \cite[Thm. 4.2.9]{KS}. We define the \textbf{quadratic Artin conductor} 
\begin{align}
Art(K)\in GW(k)[1/2]
\end{align}
to be the element defined by the pair 
\begin{align}
\left(-\int_ZC_X^Z(K), Art(K_{\textrm{et}})\right)\in GW(k)[1/2p]\times\mathbb{Z}
\end{align}
by virtue of the square~\eqref{eq:GWp}.
\end{definition}

\subsection{}
The quadratic Artin conductor is well-defined even in characteristic $2$. In characteristic different from $2$, if $Z$ has everywhere codimension at least $2$ in $X$, or if $X$ is odd-dimensional, then we have $Art(K)\in GW(k)$.

\subsection{}
If $X$ is odd-dimensional, then $C_X(\mathbbold{1}_X,\mathbf{H}_{\rm{W}}\Lambda)=e(T_{X/k},\mathbf{H}_{\rm{W}}\Lambda)=0$ . In this case we circumvent the definition of rank in~\ref{num:rkwd} and define directly
\begin{align}
\label{eq:oddrank}
\operatorname{rk}(K,\mathbf{H}_{\rm{MW}}\Lambda)\cdot C_X(\mathbbold{1}_X,\mathbf{H}_{\rm{MW}}\Lambda)
=
\operatorname{rk}(K_{\textrm{et}})\cdot \phi^*e(T_{X/k},\mathbf{H}\Lambda)
\in
\mathbf{H}_{\rm{MW}}\Lambda_0(X/k).
\end{align}
Here $\operatorname{rk}(K_{\textrm{et}})\in\mathbb{Z}$ is given by the rank of the \'etale realization of $K$, and $\phi^*:\mathbf{H}\Lambda_0(X/k)\to\mathbf{H}_{\rm{MW}}\Lambda_0(X/k)$ is induced by the canonical map $\phi:\mathbf{H}_{\rm{MW}}\Lambda\to\mathbf{H}\Lambda$ which sends $\eta$ to $0$. In other words, the right-hand side of~\eqref{eq:oddrank} only takes into account the Chow part of the Euler class.

\subsection{}
\label{num:goss}
Let $p:X\to \operatorname{Spec}(k)$ be a smooth and proper morphism with $X$ connected, and let $Z$ be a nowhere dense closed subscheme of $X$ with open complement $U$. Let $K\in\mathbf{SH}_c(X)$ be such that $K_{|U}$ is dualizable. We define $\chi(p_*K)$ and $p_*(\operatorname{rk}(K)\cdot e(T_{X/k}))$ in $GW(k)[1/2]$ in a way similar to Definition~\ref{def:quadSwan}, for the other two terms that appear in~\eqref{eq:pfgos}. In the following cases these elements are defined in $GW(k)$ without inverting $2$:
\begin{enumerate}
\item
If $k$ has characteristic different from $2$, then $\chi(p_*K)\in GW(k)$.
\item
If $k$ has characteristic different from $2$ and $Z$ has everywhere codimension at least $2$ in $X$, then $p_*(\operatorname{rk}(K)\cdot e(T_{X/k}))\in GW(k)$.
\end{enumerate}

\begin{corollary}
\label{cor:GOS2}
In the setting of~\ref{num:goss}, the following equality holds:
\begin{align}
\label{eq:GOS2}
\chi(p_*K)
=
p_*(\operatorname{rk}(K)\cdot e(T_{X/k}))-Art(K)\in GW(k)[1/2].
\end{align}
\end{corollary}

\subsection{}
We also have the following result without inverting $2$:
\begin{corollary}
\label{cor:GOS3}
In the setting of~\ref{num:goss}, assume that $k$ has characteristic different from $2$, and $X$ is odd-dimensional. Then the following equality holds:
\begin{align}
\label{eq:GOS3}
\chi(p_*K)
=
\operatorname{rk}(K_{\textrm{et}})\cdot \chi(X/k)-Art(K)\in GW(k).
\end{align}
\end{corollary}
\proof
The equality~\eqref{eq:GOS3} holds in $\mathbb{Z}$ by taking \'etale realization, and it suffices to prove the formula in $W(k)$. Since $X$ is odd-dimensional, the Witt-valued Euler class $e(T_{X/k},\mathbf{H}_{\rm{W}}\Lambda)$ is trivial and $\chi(X/k)\in GW(k)$ is hyperbolic (\cite[Cor. 4.2]{LevA}). In this case the equality $\chi(p_*K)=-Art(K)$ follows from taking degree on both sides of~\eqref{eq:CZXW}.
\endproof

\subsection{}
The formulas~\eqref{eq:GOS2} and~\eqref{eq:GOS3} give quadratic refinements of the classical Grothendieck-Ogg-Shafarevich formula (\cite[X Thm. 7.1]{SGA5}). Note that they are unconditional and do not require resolution of singularities as stated in condition~\ref{num:RS}.

\section{The non-acyclicity class}

\label{num:NAclass}

\subsection{}
In general, the localized characteristic class in Definition~\ref{def:swan} fails to commute with proper push-forwards. In order to study its functorialities we introduce a generalization to the relative situation. First recall the notion of \emph{universal local acyclicity} in motivic homotopy:
\begin{definition}[\textrm{\cite[Def. 2.1.7]{JY}}]
\label{def:locacy}
Let $f:X\to S$ be a morphism of schemes and $K\in \mathcal{SH}(X)$. We say that $K$ is \textbf{strongly locally acyclic} over $S$ if for any Cartesian square
\begin{align}
\begin{split}
  \xymatrix@=10pt{
    Y \ar[r]^-{q} \ar[d]_-{g} & X \ar[d]^-{f}\\
    T \ar[r]_-{p} & S
  }
\end{split}
\end{align}
and any object $L\in\mathcal{SH}(T)$, the canonical map $K\otimes f^*p_*L\to q_*(q^*K\otimes g^*L)$ is an isomorphism. We say that $K$ is \textbf{universally strongly locally acyclic} (abbreviated as \textbf{USLA}) over $S$ if for any morphism $T\to S$, the base change $K_{|X\times_ST}$ is strongly locally acyclic over $T$.
\end{definition}

\subsection{}
The following properties hold (see \cite{JY}, \cite{LZ}, \cite[\S3.2]{Pre}):
\begin{lemma}
\label{lm:latensd}
If $K$ is USLA over $S$, then for any dualizable object $A$, $K\otimes A$ is USLA over $S$.
\end{lemma}

\begin{lemma}
\label{lm:LZ218}
If $K\in\mathbf{SH}(X)$ is USLA over $S$, then $\mathbb{D}_{X/S}(K)$ is USLA over $S$.
\end{lemma}

\begin{lemma}
\label{lm:laprod}
Let $S$ be a scheme and let $f_1:X_1\to Y_1$, $f_2:X_2\to Y_2$ be two morphisms of $S$-schemes. For $i=1,2$, let $K_i\in\mathbf{SH}(X_i)$ be an object which is USLA over $Y_i$. Then the object $K_1\boxtimes_SK_2$ is USLA over $Y_1\times_SY_2$.
\end{lemma}
\proof
Consider the following commutative diagram 
\begin{align}
\begin{split}
  \xymatrix{
   X_1 \ar[d]_-{f_1} \ar@{}[rd]|{\Delta_1} & X_1\times_S X_2 \ar[l]_-{p_1} \ar[rd]^-{p_2} \ar[d]^-{f_1\times id} & \\
   Y_1 & Y_1\times_S X_2 \ar[l]_-{\underline{p'_1}} \ar[r]^-{\underline{p_2}} \ar[d]_-{id\times f_2} \ar@{}[rd]|{\Delta_2} & X_2 \ar[d]^-{f_2}\\
   & Y_1\times_S Y_2 \ar[lu]^-{p_1'} \ar[r]^-{p_2'} & Y_2
  }
\end{split}
\end{align}
where both squares are Cartesian. Let $V\xrightarrow{p}Y_1\times_S Y_2$ be a morphism and form the following Cartesian square
\begin{align}
\begin{split}
  \xymatrix@=10pt{
  W \ar[d]_-{g_1} \ar[r]^-{r} & X_1\times_S X_2 \ar[d]^-{f_1\times id} \\
  V_1 \ar[d]_-{g_2} \ar[r]^-{q} & Y_1\times_S X_2 \ar[d]^-{id\times f_2} \\
  V \ar[r]^-{p} & Y_1\times_S Y_2.
  }
\end{split}
\end{align}
Then for any object $L\in\mathbf{SH}(V)$, our assumptions imply the following isomorphism
\begin{align}
\begin{split}
&(K_1\boxtimes_SK_2)\otimes(f_1\times id)^*(id\times f_2)^*p_*L\\
=
&p_1^*K_1\otimes(f_1\times id)^*(\underline{p_2}^*K_2\otimes(id\times f_2)^*p_*L)\\
\simeq
&p_1^*K_1\otimes(f_1\times id)^*q_*(q^*\underline{p_2}^*K_2\otimes g_2^*L)\\
\simeq
&r_*(r^*p_1^*K_1\otimes g_1^*(q^*\underline{p_2}^*K_2\otimes g_2^*L))\\
=
&r_*(r^*(K_1\boxtimes_SK_2)\otimes g_1^*g_2^*L)
\end{split}
\end{align}
which shows that $K_1\boxtimes_SK_2$ is SLA over $Y_1\times_SY_2$. The same property also holds after any base change, so $K_1\boxtimes_SK_2$ is USLA over $Y_1\times_SY_2$.
\endproof

\subsection{}
For any Cartesian square of schemes
\begin{align}
\label{eq:diagDelta}
\begin{split}
  \xymatrix@=10pt{
   Y \ar[r]^-{q} \ar[d]_-{g} \ar@{}[rd]|{\Gamma} & X \ar[d]^-{f} \\
    T \ar[r]_-{p} & S
    }
  \end{split}
\end{align}
where $p$ is separated of finite type, and any objects $K\in\mathbf{SH}(X)$ and $L\in\mathbf{SH}(T)$, there is a canonical natural transformation
\begin{align}
\label{eq:YZ213}
q^*K\otimes g^*p^!L\to q^*K\otimes q^!f^*L\to q^!(K\otimes f^*L).
\end{align}

\begin{lemma}
\label{lm:Ill25}
If $K$ is USLA over $S$, then the map~\eqref{eq:YZ213} is an isomorphism. 
\end{lemma}
\proof
The problem is Zariski local in $T$, so we may assume that $p$ is quasi-projective; the problem is also stable under compositions of the morphism $p$, so we may further assume that $p$ is either smooth or a closed immersion. If $p$ is smooth, then the map~\eqref{eq:YZ213} is an isomorphism by smooth base change, and it remains to show the case where $p$ is a closed immersion.

Assume that $p$ is a closed immersion, and denote by $j$ and $k$ the complementary open immersions of $p$ and $q$. Since $K$ is USLA over $S$, there is a canonical isomorphism
\begin{align}
\label{eq:jUSLA}
K\otimes f^*j_*j^*L\simeq k_*k^*(K\otimes f^*L).
\end{align}
We have the following commutative diagram
\begin{align}
\begin{split}
  \xymatrix@=10pt{
   K\otimes f^*p_*p^!L \ar[r]^-{\sim} \ar[d]_-{} & K\otimes q_*g^*p^!L \ar[r]^-{\sim} & q_*(q^*K\otimes g^*p^!L) \ar[r]^-{\eqref{eq:YZ213}} & q_*q^!(K\otimes f^*L) \ar[d]_-{} \\
    K\otimes f^*L \ar@{=}[rrr]^-{} \ar[d]_-{} & & & K\otimes f^*L \ar[d]_-{} \\
    K\otimes f^*j_*j^*L \ar[rrr]^-{\eqref{eq:jUSLA}}_-{\sim} & & & k_*k^*(K\otimes f^*L) \\
    }
  \end{split}
\end{align}
where both columns are part of the localization cofiber sequence. By the five lemma, since the functor $q_*$ is conservative, we deduce that the map~\eqref{eq:YZ213} is an isomorphism when $p$ is a closed immersion, which finishes the proof.
\endproof

\begin{definition}
Given a diagram of the form~\eqref{eq:diagDelta}, we denote by $\Gamma^{\Delta}:\mathbf{SH}(X)\to\mathbf{SH}(Y)$ the cofiber of the map
\begin{align}
q^*(-)\otimes g^*p^!\mathbbold{1}_S\xrightarrow{\eqref{eq:YZ213}}q^!(-).
\end{align}
In particular, Lemma~\ref{lm:Ill25} implies that if $K$ is USLA over $S$ then $\Gamma^{\Delta}(K)=0$.
\end{definition}

\begin{remark}
If $f$ is the identity morphism, then $\Gamma^{\Delta}=p^{\Delta}$ is the map defined in~\ref{num:JY614}.
\end{remark}


 

\subsection{}
\label{num:laset}
Let $X\xrightarrow{f}Y\xrightarrow{g}S$ be a diagram of schemes with $g$ smooth, and consider the following Cartesian diagram
\begin{align}
\label{eq:YZ252}
\begin{split}
  \xymatrix@=10pt{
     X\ar@{=}[r]\ar[d]_-{\delta_{X/Y}} & X\ar[d]^-{\delta_{X/S}} \\
   X\times_YX \ar[r]^-{\delta_{X/Y/S}} \ar[d]_-{} \ar@{}[rd]|{\Gamma_{X/Y/S}} & X\times_SX \ar[d]^-{f\times_Sf} \\
    Y \ar[r]_-{\delta_{Y/S}} & Y\times_SY.
    }
  \end{split}
\end{align}

\begin{corollary}
\label{cor:YZ253}
Assume the setting of~\ref{num:laset}.
Let $K,L\in\mathbf{SH}_c(X)$ be objects which are USLA over $Y$.  Then we have 
\begin{align}
\Gamma_{X/Y/S}^\Delta(\mathbb{D}_{X/S}(K)\boxtimes_SL)=0.
\end{align}
In other words, there is a canonical isomorphism
\begin{align}
\mathbb{D}_{X/Y}(K)\boxtimes_YL
\simeq
\delta_{X/Y/S}^!(\mathbb{D}_{X/S}(K)\boxtimes_SL).
\end{align}
\end{corollary}
\proof
By Lemma~\ref{lm:LZ218}, the object $\mathbb{D}_{X/Y}(K)$ is USLA over $Y$. Since $Y$ is smooth over $S$, by Lemma~\ref{lm:latensd}, the object $\mathbb{D}_{X/S}(K)$ is USLA over $Y$. By Lemma~\ref{lm:laprod}, the object $\mathbb{D}_{X/S}(K)\boxtimes_SL$ is USLA relatively to $f\times_Sf:X\times_SX\to Y\times_SY$. We conclude by applying Lemma~\ref{lm:Ill25}.
\endproof

\subsection{}
Assume the setting of~\ref{num:laset}. Since the object $\delta_{Y/S}^!\mathbbold{1}_{Y\times_SY}$ is dualizable (in fact $\otimes$-invertible), the following canonical map is an isomorphism:
\begin{align}
f^!(-)\otimes f^*\delta_{Y/S}^!\mathbbold{1}_{Y\times_SY}
\to
f^!(-\otimes\delta_{Y/S}^!\mathbbold{1}_{Y\times_SY})
\end{align}
from which we deduce an isomorphism of functors
\begin{align}
\label{eq:Kalt}
\delta_{X/Y}^*\Gamma_{X/Y/S}^\Delta\delta_{X/S,*}f^!
\simeq
f^!\delta_{Y/S}^\Delta\delta_{Y/S,*}.
\end{align}

\begin{definition}
Assume the setting of~\ref{num:laset}.
We denote 
\begin{align}
\label{eq:Kreldef}
\mathcal{K}_{X/Y/S}
=
f^!\delta_{Y/S}^\Delta\delta_{Y/S,*}\mathcal{K}_{Y/S}
\overset{\eqref{eq:Kalt}}{\simeq}
\delta_{X/Y}^*\Gamma_{X/Y/S}^\Delta\delta_{X/S,*}\mathcal{K}_{X/S}
\in\mathbf{SH}_c(X).
\end{align}
\end{definition}

\subsection{}
There is a cofiber sequence
\begin{align}
\mathcal{K}_{X/Y}
\to
\mathcal{K}_{X/S}
\to
\mathcal{K}_{X/Y/S}.
\end{align}

\subsection{}
If $g:W\to X$ is a morphism, by~\eqref{eq:Kreldef} we have
\begin{align}
\label{eq:Krelcomp}
\mathcal{K}_{W/Y/S}=g^!\mathcal{K}_{X/Y/S}.
\end{align}

\subsection{}
\label{num:larset}
Assume the setting of~\ref{num:laset}.
Let $i:Z\to X$ be a closed immersion. Let $K\in\mathbf{SH}_c(X)$ be an object which is USLA over $S$ and such that $K_{|X-Z}$ is USLA over $Y$. 
By Corollary~\ref{cor:YZ253}, the following canonical map is an isomorphism:
\begin{align}
\label{eq:lau}
\delta_{X/Y}^*\Gamma_{X/Y/S}^\Delta(\mathbb{D}_{X/S}(K)\boxtimes_SK)
\simeq
i_*i^*\delta_{X/Y}^*\Gamma_{X/Y/S}^\Delta(\mathbb{D}_{X/S}(K)\boxtimes_SK)
\end{align}
\begin{definition}
\label{def:na}
In the setting of~\ref{num:larset}, we define a class
\begin{align}
B^{Z}_{X/Y/S}(K)
\in 
[i_*\mathbbold{1}_Z,\mathcal{K}_{X/Y/S}]
\simeq
[\mathbbold{1}_Z,\mathcal{K}_{Z/Y/S}]
\end{align}
as the composition
\begin{align}
\begin{split}
i_*\mathbbold{1}_Z
&\to
i_*i^*\delta_{X/Y}^*\Gamma_{X/Y/S}^\Delta\delta_{X/S,*}\mathbbold{1}_X
\to
i_*i^*\delta_{X/Y}^*\Gamma_{X/Y/S}^\Delta(\mathbb{D}_{X/S}(K)\boxtimes_SK)\\
&\overset{\eqref{eq:lau}}{\simeq}
\delta_{X/Y}^*\Gamma_{X/Y/S}^\Delta(\mathbb{D}_{X/S}(K)\boxtimes_SK)
\simeq
\delta_{X/Y}^*\Gamma_{X/Y/S}^\Delta(K\boxtimes_S\mathbb{D}_{X/S}(K))\\
&\to
\delta_{X/Y}^*\Gamma_{X/Y/S}^\Delta\delta_{X/S,*}\mathcal{K}_{X/S}
\overset{\eqref{eq:Kreldef}}{\simeq}
\mathcal{K}_{X/Y/S}.
\end{split}
\end{align}
The class $B_X^Z(K)$ in~\eqref{eq:mapZ} is nothing but the class $B_{X/X/k}^Z(K)$.
\end{definition}

\begin{proposition}[Additivity]
If $L\to M\to N$ is a distinguished triangle in $\mathbf{SH}_c(X)$, then one has
\begin{align}
B_{X/Y/S}^Z(M)=B_{X/Y/S}^Z(L)+B_{X/Y/S}^Z(N).
\end{align}
\end{proposition}
The proof is analogous to Lemma~\ref{lm:addB}.

\subsection{}
Consider a commutative diagram
\begin{align}
\label{eq:diagGammaCart}
\begin{split}
  \xymatrix@=10pt{
   W \ar[r]^-{} \ar[d]_-{p} & U \ar[r]^-{} \ar[d]_-{q} & T \ar[d]^-{r} \\
    X \ar[r]_-{f} & Y \ar[r]_-{g} & S
    }
  \end{split}
\end{align}
where the square on the right is Cartesian.
One can associate a Cartesian square
\begin{align}
\begin{split}
  \xymatrix@=10pt{
    W\times_UW \ar[r]^-{} \ar[d]_-{p\times_Yp} & W\times_TW \ar[d]^-{p\times_Sp} \\
    X\times_YX \ar[r]_-{} & X\times_SX
    }
  \end{split}
\end{align}
which is obtained by base change from the Cartesian square
\begin{align}
\begin{split}
  \xymatrix@=10pt{
    U \ar[r]^-{\delta_{U/T}} \ar[d]_-{q} & U\times_TU \ar[d]^-{q\times_Sq} \\
    Y \ar[r]_-{\delta_{Y/S}} & Y\times_SY.
    }
  \end{split}
\end{align}
We have a natural transformation
\begin{align}
q^*\delta_{Y/S}^!\to \delta_{U/T}^!(q\times_Sq)^*
\end{align}
from which we deduce a natural transformation
\begin{align}
\label{eq:natGamma}
(p\times_Yp)^*\Gamma_{X/Y/S}^\Delta
\to
\Gamma_{W/U/T}^\Delta(p\times_Sp)^*.
\end{align}
If the left square in~\eqref{eq:diagGammaCart} is also Cartesian, we obtain a map
\begin{align}
\label{eq:natrelK}
p^*\mathcal{K}_{X/Y/S}
\to
\mathcal{K}_{W/U/T}.
\end{align}

\begin{proposition}[Base change]
Assume the setting of Definition~\ref{def:na}. 
Consider a Cartesian square
\begin{align}
\begin{split}
  \xymatrix@=10pt{
   Z_T \ar[r]^-{i_T} \ar[d]_-{p_Z}  & X_T \ar[r]^-{f_T} \ar[d]_-{p_X} & Y_T \ar[r]^-{g_T} \ar[d]_-{p_Y} & T \ar[d]^-{p} \\
   Z \ar[r]_-{i} & X \ar[r]_-{f} & Y \ar[r]_-{g} & S.
    }
  \end{split}
\end{align}
Then the following diagram is commutative:
\begin{align}
\begin{split}
  \xymatrix{
     p_X^*i_*\mathbbold{1}_Z \ar[r]^-{p_X^*B^{Z}_{X/Y/S}(K)} \ar[d]_-{\wr} & p_X^*\mathcal{K}_{X/Y/S} \ar[d]^-{\eqref{eq:natrelK}} \\
    i_{T*}\mathbbold{1}_{Z_T} \ar[r]^-{B^{Z_T}_{X_T/Y_T/T}(p_X^*K)} & \mathcal{K}_{X_T/Y_T/T}
    }
  \end{split}
\end{align}
\end{proposition}
\proof
Use the natural transformation
\begin{align}
(p_X\times_Yp_X)^*\Gamma_{X/Y/S}^\Delta
\to
\Gamma_{X_T/Y_T/T}^\Delta(p_X\times_Sp_X)^*
\end{align}
induced by~\eqref{eq:natGamma}.
\endproof

\begin{proposition}[Push-forward]
Assume the setting of Definition~\ref{def:na}. 
Consider a commutative diagram
\begin{align}
\begin{split}
  \xymatrix@=10pt{
   Z \ar[r]^-{i} \ar[d]_-{r} & X \ar[rd]^-{f} \ar[d]_-{p} & &  \\
   Z' \ar[r]_-{i'} & X' \ar[r]_-{q} & Y \ar[r]_-{g} & S.
    }
  \end{split}
\end{align}
where $p$ is proper and both $i$ and $i'$ are closed immersions.
Then the following diagram is commutative:
\begin{align}
\begin{split}
  \xymatrix{
     i'_{*}\mathbbold{1}_{Z'} \ar[r]^-{B^{Z'}_{X'/Y/S}(p_*K)} \ar[d]_-{} & \mathcal{K}_{X'/Y/S} & p_*p^!\mathcal{K}_{X'/Y/S} \ar[l]^-{} \\
    i'_{*}r_*\mathbbold{1}_{Z} \ar@{=}[r]_-{} & p_{*}i_*\mathbbold{1}_{Z} \ar[r]^-{p_*B^{Z}_{X/Y/S}(K)} & p_*\mathcal{K}_{X/Y/S} \ar@{=}[u]_-{\eqref{eq:Krelcomp}}
    }
  \end{split}
\end{align}

\end{proposition}
\proof
Use the natural transformation
\begin{align}
(p\times_Yp)^*\Gamma_{X'/Y/S}^\Delta
\to
\Gamma_{X/Y/S}^\Delta(p\times_Sp)^*
\end{align}
induced by~\eqref{eq:natGamma}.
\endproof

\begin{proposition}[Pull-back]
Assume the setting of Definition~\ref{def:na}. 
Consider a commutative diagram
\begin{align}
\begin{split}
  \xymatrix@=10pt{
   V \ar[r]^-{k} \ar[d]_-{r} & W \ar[rd]^-{q} \ar[d]_-{p} & &  \\
   Z \ar[r]_-{i} & X \ar[r]_-{f} & Y \ar[r]_-{g} & S.
    }
  \end{split}
\end{align}
where $p$ is \'etale and both $i$ and $k$ are closed immersions.
Then the following diagram is commutative:
\begin{align}
\begin{split}
  \xymatrix{
     p^*i_*\mathbbold{1}_Z \ar[r]^-{p^*B^{Z}_{X/Y/S}(K)} \ar[d]_-{} & p^*\mathcal{K}_{X/Y/S} \ar@{=}[d]^-{\eqref{eq:Krelcomp}} \\
    k_*\mathbbold{1}_{V} \ar[r]^-{B^{V}_{W/Y/S}(p^*K)} & \mathcal{K}_{W/Y/S}.
    }
  \end{split}
\end{align}


\end{proposition}
\proof
Use the natural transformation
\begin{align}
(p\times_Yp)^*\Gamma_{X/Y/S}^\Delta
\to
\Gamma_{W/Y/S}^\Delta(p\times_Sp)^*
\end{align}
induced by~\eqref{eq:natGamma}.
\endproof

\subsection{}
If the morphism $Y\to S$ factors as $Y\xrightarrow{k}T\xrightarrow{h}S$ with both $k$ and $h$ smooth, then there is a canonical map
\begin{align}
\label{eq:Krelfac}
\mathcal{K}_{X/Y/S}
\to
\mathcal{K}_{X/T/S}.
\end{align}
Indeed, this map is constructed from the natural transformation 
\begin{align}
\delta_{Y/S}^\Delta\delta_{Y/S,*}k^!
\to
k^!\delta_{T/S}^\Delta\delta_{T/S,*}
\end{align}
which is deduced from the commutative square
\begin{align}
\label{eq:YZ252}
\begin{split}
  \xymatrix@=10pt{
     Y \ar[r]^-{\delta_{Y/T}} \ar[rd]_-{k} & Y\times_TY \ar[r]^-{} \ar[d]^-{} \ar@{}[rd]|{\Gamma_{Y/T/S}} & Y\times_SY \ar[d]^-{k\times_Sk} \\
    & T \ar[r]_-{\delta_{T/S}} & T\times_ST.
    }
  \end{split}
\end{align}
The map~\eqref{eq:Krelfac} is an isomorphism when $k$ is \'etale. 

\begin{lemma}
\label{lm:smlocb}
If the morphism $Y\to S$ factors as $Y\xrightarrow{k}T\xrightarrow{h}S$ with both $k$ and $h$ smooth, then 
the following diagram is commutative:
\begin{align}
\begin{split}
  \xymatrix{
     i_*\mathbbold{1}_Z \ar[r]^-{B^{Z}_{X/Y/S}(K)} \ar[rd]_-{B^{Z}_{X/T/S}(K)} & \mathcal{K}_{X/Y/S} \ar[d]^-{\eqref{eq:Krelfac}} \\
    & \mathcal{K}_{X/T/S}
    }
  \end{split}
\end{align}
\end{lemma}
\proof
This follows from the existence of a natural transformation 
\begin{align}
\Gamma_{X/Y/S}^\Delta
\to
\delta_1^*\Gamma_{X/T/S}^\Delta
\end{align}
where $\delta_1$ is the morphism $X\times_YX\to X\times_TX$.
\endproof

\begin{lemma}
\label{lm:smlocb}
If the morphism $Y\to S$ factors as $Y\xrightarrow{k}T\xrightarrow{h}S$ with $k$ smooth and $h$ the composition of an \'etale morphism with a universal homeomorphism, then 
the following diagram is commutative:
\begin{align}
\begin{split}
  \xymatrix{
     i_*\mathbbold{1}_Z \ar[r]^-{B^{Z}_{X/Y/T}(K)} \ar[rd]_-{B^{Z}_{X/Y/S}(K)} & \mathcal{K}_{X/Y/T} \ar[d]^-{\wr} \\
    & \mathcal{K}_{X/Y/S}.
    }
  \end{split}
\end{align}
\end{lemma}



\end{document}